\documentclass{elsart}
\usepackage{natbib}
\usepackage{latexsym}
\usepackage[centertags]{amsmath}
\usepackage{amsfonts}
\usepackage{amssymb}

\newtheorem{exa}[thm]{Example}

\makeatletter\@addtoreset{equation}{section}\makeatother

\def\QQ{\mathbb{Q}}

\def\CC{\mathbb{C}}

\def\poly{\mbox{\sf Poly}}
\def\frr{\mbox{\sf Frac}}

\newcommand{\ceiling}[1]{\,\lceil #1 \rceil\,}

\begin{document}

\begin{frontmatter}

\title{A Fast Algorithm for Partial Fraction
Decompositions}

\author{Guoce Xin}
\ead{maxima@brandeis.edu}

\address{Department
of Mathematics\\
Brandeis University\\
Waltham MA 02454-9110}

\begin{abstract}
We obtain two new algorithms for partial fraction decompositions;
the first is over algebraically closed fields, and the second is
over general fields. These algorithms takes $O(M^2)$ time, where
$M$ is the degree of the denominator of the rational function. The
new algorithms use less storage space, and are suitable for
parallel programming. We also discuss full partial fraction
decompositions.
\end{abstract}


\end{frontmatter}

\newcommand\xmod[1]{\langle #1 \rangle}

\section{Introduction}

The partial fraction decomposition of a one-variable rational
function is very useful in mathematics. For example, it is crucial
to obtain the partial fraction decomposition of a rational
function in order to integrate it. Kovacic's algorithm
\cite{kovacic} for solving the differential equation
$y''(x)+r(x)y(x)=0$, where $r(x)$ is rational, requires the full
partial fraction expansion of $r(x)$ over the complex numbers.

The classical algorithm for partial fraction expansion relies on
the following theorem. To make it simple, we consider rational
functions in $\CC(t)$.

\begin{thm}
If $a_1,\dots ,a_n$ are $n$ distinct complex numbers, $m_1,\dots
,m_n$ are positive integers, and the degree of $p(t)$ is less than
$m_1+\cdots +m_n$, then there are unique complex numbers
$A_{i,j}$, where $1\le i\le n$ and $1\le j\le m_i$, such that
\begin{equation}
\label{e-parfrac1} \frac{p(t)}{(t-a_1)^{m_1}\cdots (t-a_n)^{m_n}}
= \sum_{i=1}^n  \sum_{j=1}^{m_i} \frac{A_{i,j}}{(t-a_i)^j}.
\end{equation}
\end{thm}

The classical algorithm multiplies both sides by the denominator,
and then equates coefficients to solve a large system of linear
equations for the $A_{i,j}$'s.


The key observation for our first algorithm is that linear
transformations will preserve the structure of a partial fraction
expansion. This reduces the problem of finding all $A_{i,j}$ to
finding only $A_{1,j}$ and assuming that $a_1=0$. For finding
$A_{1,j}$, we use the unique Laurent series expansion at $t=0$.

Denote by $F(t)$ the left side of equation \eqref{e-parfrac1}. Let
$M$ be the degree of the denominator of $F(t)$, which is
$m_1+m_2+\cdots +m_k$. Compared with the classical algorithm for
obtaining the partial fraction decomposition of $F(t)$, our new
algorithm has three advantages. This comparison is under the
assumption of fast multiplication of (usually rational) numbers.
In the following, when we say that an algorithm takes $O(M)$ time,
we mean that the algorithm will do
  $O(M)$ multiplications.

\begin{enumerate}
\item The new algorithm is fast. The classical algorithm needs to
solve $M$ linear equations in $M$ unknowns, which takes
$\Omega(M^3)$ time using the Gaussian elimination algorithm. See
\cite[p. 540, Property 37.1]{sedgewick}. But our algorithm only
takes $O(M^2)$ time.

\item The new algorithm requires little storage space.
The classical algorithm needs to record all of the $M^2$
coefficients in these $M$ linear equations. So the storage space
is about $\Omega(M^2)$. But our new algorithm needs only to record
two polynomials of degree $m$, where $m$ is the maximum of the
$m_i$'s. So the storage space is at most $O(M)$.

\item The new algorithm computes the partial fraction
expansion at different $a_i$'s separately, so it is more suitable
for parallel programming.
\end{enumerate}


For partial fraction decompositions in a general field $K$, we
also have a fast algorithm by working in some quotient rings. The
new algorithm applies to finding the full partial fraction
decompositions and to evaluating generalized Dedekind sums. The
theory gives an efficient algorithm for MacMahon's partition
analysis \cite[Ch. 2.5]{xinthesis}.

\section{Partial Fraction Decompositions Over a General Field}

Let $K$ be any field, and let $t$ be a variable. It is well-known
that the ring of polynomials $K[t]$ has many nice properties. Here
we use the fact that $K[t]$ is a unique factorization domain. We
will use a quotient ring to derive a formula for partial fraction
decomposition, which is the basis of our new algorithms.

In what follows, the degree of an element $r\in K[t]$, denoted by
$\deg(r)$, is the degree of $r$ as a polynomial in $t$. The degree
of the $0$ polynomial is treated as $-\infty$. We start with the
division theorem in $K[t]$.

\begin{thm}[Division Theorem]\label{t-division}
Let $D ,N \in K[t]$ and suppose $D\ne 0$. There is a unique pair
$(p,r)$ such that $p,r\in K[t]$, $N=Dp+r$, and $\deg r< \deg D$.
\end{thm}
The $r$ above is called the {\em remainder} of $N$ when divided by
$D$.

A rational function $N/D$ with $N,D\in K[t]$ is said to be {\em
proper}
 if $\deg N <\deg D$. A proper rational function is
simply called a {\em proper fraction}. The unit $1$ is not proper,
but $0$ is considered to be proper. It is clear that the sum of
proper fractions is  a proper fraction, and the product of proper
fractions is  a proper fraction. But the set of all proper
fractions does not form a ring, for $1$ does not belong to it.

By the division theorem, any rational function $N/D$ can be
uniquely written as the sum of a polynomial and a proper fraction.
Such a decomposition is called a {\em ppfraction} (short for
polynomial and proper fraction) of $N/D$. If $N=Dp+r$ with
$\deg(r)<\deg(D)$, then $N/D=p+r/D$ is a ppfraction. We denote by
$\poly(N/D)$ the polynomial part of $N/D$, and by $\frr(N/D)$ the
fractional part of $N/D$.

Recall the following well-known result in algebra.
\begin{lem}\label{l-ppfraction-k}
Let $N,D\in K[t]$ with $D\ne 0$. If $D=D_1\cdots D_k$ is a
factorization of $D$ in $K[t]$, and all the $D_i$ are pairwise
relatively prime, then $N/D$ can be uniquely written as
\begin{equation}\label{e-ppfraction-k}
\frac{N}{D}= p+\frac{r_1}{D_1}+\cdots +\frac{r_k}{D_k},
\end{equation}
where $r_i$ is a polynomial of degree smaller than $\deg (D_i)$
for all $i$, and $p$ equals the polynomial part of $N/D$.
\end{lem}
We call the above decomposition the \emph{ppfraction expansion} of
$N/D$ \emph{with respect to} $(D_1,\dots ,D_k)$. To find such
decomposition, it suffices to find $p$ and $r_1,\dots ,r_k$. It is
easy to find $p$. In finding $r_i$, it is convenient to use the
concept of quotient ring $K[t]/\langle D_i\rangle $, where
$\langle D_i\rangle $ is the ideal generated by $D_i$.

Recall that $D'+\langle D\rangle $ has a multiplicative inverse in
$K[t]/\langle D\rangle $ if and only if $D'$ is relatively prime
to $D$. Moreover, if $\frac{1}{DD'}=\frac{r}{D}+\frac{r'}{D'}$,
which is not necessarily a ppfraction expansion, then
$1/D'+\langle D\rangle =r+\langle D\rangle $. The last fact we
will use is that for any polynomials $N_1$ and $N_2$, $N_1+\langle
D\rangle =N_2+\langle D\rangle $ if and only if
$\frr(N_1/D)=\frr(N_2/D)$.

Suppose that $D=D_1D'$ and that $D_1$ and $D'$ are relatively
prime. Then we have a ppfraction of $N/D$ with respect to
$(D_1,D')$:
$$N/D=\poly(N/D)+r_1/D_1+r'/D'.$$
In such a decomposition, we  call $r_1/D_1$  the {\em fractional
part of $N/D$ with respect to $D_1$}, and denote it by
$\frr(N/D,D_1)$. If $D_1=(t-a)^m$ for some $a\in K$, then we
simply denote it by $\frr(N/D,t=a)$. Clearly $\frr(N/D,1)$ is
always $0$, and $\frr(N/D,D_1)$ is always a proper fraction with
denominator $D_1$. Also we have the following relation:
$$r_1+\langle D_1\rangle =\frac{N}{D'}+\langle D_1\rangle .$$
Thus to find $r_1$, we pick a representative of $N/D'+\langle
D_1\rangle $, and then find its remainder when divided by $D_1$.

\begin{thm}\label{t-ppfraction-main1}
For any $N,D\in K[t]$ with $D\ne 0$, if $D_1,\dots D_k \in K[t]$
are pairwise relatively prime, and $D=D_1\cdots D_k$, then
$$\frac{N}{D}= \poly\left(\frac{N}{D}\right) +\frr \left(\frac{N}{D}, D_1 \right)+\cdots +
\frr \left(\frac{N}{D}, D_k \right)$$ is the ppfraction expansion
of $N/D$ with respect to $(D_1,\dots ,D_k)$. Moreover, if
$1/(D_1D_i)= s_i/D_1+p_i/D_i$,  then
$$\frr(N/D,D_1)
=\frr(Ns_2s_3\cdots s_k/D_1).$$
\end{thm}
\begin{pf}
For the first part, suppose that
\begin{align} \label{e-thmproof1}\frac{N}{D}=
p+\frac{r_1}{D_1}+\cdots +\frac{r_k}{D_k}
\end{align} is the ppfraction
expansion of $N/D$ with respect to $(D_1,\dots ,D_k)$. Let
$D'=D_2\cdots D_k$. Then $D_1$ and $D'$ are relatively prime and
$r_2/D_2+\cdots +r_k/D_k=r'/D'$ is a proper fraction with
denominator $D'$. By the uniqueness of ppfraction of $N/D$ with
respect to $(D_1,D')$, we have $r_1/D_1= \frr(N/D,D_1)$. Similarly
$r_i/D_i= \frr(N/D,D_i)$ for all $i$.

For the second part, multiplying both sides of equation
\eqref{e-thmproof1} by $D$, and thinking of this as an identity in
the quotient ring $K[t]/\langle D_1\rangle $, we get
$$N+\langle D_1\rangle =pD_1D'+r_1D_2\cdots D_k+\cdots +r_kD_1\cdots D_{k-1}+\langle D_1\rangle =r_1D'+\langle D_1\rangle .$$
Now multiply both sides of the above equation by $1/D'+\langle
D_1\rangle $, we get
$$r_1+\langle D_1\rangle =\frac{N}{D_2\cdots D_k}+\langle D_1\rangle =Ns_2s_3\cdots s_k+\langle D_1\rangle .$$ Therefore
$$\frac{r_1}{D_1}= \frr\left(\frac{r_1}{D_1}\right)=\frr\left(\frac{Ns_2s_3\cdots s_k}{D_1} \right).$$
\end{pf}

Theorem \ref{t-ppfraction-main1} is the basis of our new
algorithms. Let $M$ be the degree of the denominator of a rational
function. We will give an $O(M^2)$ algorithm for finding the
partial fraction decomposition based on the above theorem.

If $D=a p_1^{m_1} \cdots p_k^{m_k}$, where $a\in K$,
 is a factorization of $D$ into monic primes in $K[t]$,
then $p_1^{m_1}, \dots , p_k^{m_k}$ are pairwise relatively prime.
Let $D_i=p_i^{m_i}$, and let $r_i$ be a polynomial with
$\deg(r_i)<\deg(D_i)$. Then every $r_i/D_i$ can be uniquely
written in the form $\sum_{j=1}^{m_i} A_j/p_i^j$ with $\deg
(A_j)<\deg(p_i)$ for all $j$. The {\em partial fraction expansion}
of $N/D$ is the result of applying the above decomposition to the
ppfraction of $N/D$ with respect to $(D_1,\dots ,D_k)$. In this
case, we can use the following lemma to reduce the problem to
computing only the partial fraction expansion of $1/(p_ip_j)$ for
all $i\ne j$.

\begin{lem}\label{l-frac-power}
Let $p,q\in K[t]$ be relatively prime polynomials. If $r$ and $s$
are two polynomials such that $1/(pq)=r/p+s/q$, then for any
positive integers $m,n$,
\begin{equation}\label{e-ppfraction-pq}
\frac{1}{p^mq^n} =\frac{1}{p^m} \sum_{i=0}^{m-1} \binom{m+i}{i}
r^ns^{i}p^{i} +\frac{1}{q^n} \sum_{j=0}^{n-1} \binom{n+j}{j}
r^{j}s^m q^{j}.
\end{equation}
\end{lem}
\begin{pf}
Using the formula $1/(pq)=r/p+s/q$, we have
$$\frac1{p^mq^n} =\frac{1}{pq}\cdot \frac{1}{p^{m-1}q^{n-1}}
= \frac{r}{p^mq^{n-1}}+\frac{s}{p^{m-1}q^n}.$$ If we let
$A(m,n)=1/(p^mq^n)$, then the above equation is equivalent to
$$A(m,n)=rA(m,n-1)+sA(m-1,n).$$
Using this recursive relation, we can express $A(m,n)$ in terms of
$A(0,j)$ and $A(i,0)$, where $1\le j\le n$ and $1\le i\le m$.

Either using induction or a combinatorial argument, we can easily
get
$$A(m,n)= \sum_{i=0}^{m-1} \binom{m+i}{i} r^{n}s^{i}A(m-i,0)
+\sum_{j=0}^{n-1} \binom{n+j}{j} r^{j}s^m A(0,n-j).$$ Equation
\eqref{e-ppfraction-pq} is just a restatment of the above
equation.
\end{pf}

\section{Partial Fraction Decompositions in $\CC(t)$}

In this section, $K$ is an algebraically closed field (e.g., the
field of complex numbers $\CC$). Partial fraction decomposition in
this situation is simple, since every polynomial in $K[t]$ can be
written as a product of linear factors $t-a$ for $a\in K$.

The key idea to our new algorithm is that linear transformation
will not change the structure of a partial fraction decomposition.
This can be illustrated by the following example.

The partial fraction expansion of $f(t)$ is $A/(t-a)+B/(t-b)$ if
and only if the partial fraction expansion of $f(t+c)$ is
$A/(t+c-a)+B/(t+c-b)$. So we can compute the partial fraction
expansion of $f(t+a)$, and after that, replace $t$ with $t-a$.

\vspace{3mm} Let $b\in K$ and let $\tau_b$ be the transformation
defined by $\tau_b\  f(t)=f(t+b)$ for any $f(t)\in K[t]$ or
$f(t)\in K(t)$. Then $\tau_b$ is clearly an automorphism of $K[t]$
and of $K(t)$, and its inverse is $\tau_{-b}$. The following
properties can be easily checked for any $p,q\in K[t]$ and $b\in
K$.
\begin{enumerate}
\item  $p$ is prime in $K[t]$ if and only if
$\tau_b\  p$ is.
\item $\tau_b\  \gcd(p,q) =\gcd(\tau_b\  p,\tau_b\  q)$.
\item $\deg (\tau_b\ p)=\deg(p)$.
\item $p/q$ is a proper fraction if and only if $\tau_b\ p/q$ is.
\end{enumerate}
Thus for any $N,D\in K[t]$ with $D\ne 0$, $N/D= p+r_1/D_1+\cdots
+r_k/D_k$ is the ppfraction expansion of $N/D$ if and only if
$\tau_b\  N/D= (\tau_b\  p)+(\tau_b\  r_1/D_1)+\cdots +(\tau_b\
r_k/D_k)$ is a ppfraction expansion of $\tau_b\  N/D$. The partial
fraction expansion can be obtained by first computing the partial
fraction expansion of $\tau_b\  N/D$, then applying $\tau_{-b}$ to
the result. Choosing $b$ appropriately can simplify the
computation. The above argument gives us the following lemma.
\begin{lem}\label{l-ppfraction-tr}
For any $N,D,D_1\in K[t]$ with $D\ne 0$, $D/D_1\in K[t]$, and
\\
$\gcd (D_1,D/D_1)=1$, we have
$$\frr(N/D,D_1)=\tau_{-b}\  \frr(\tau_b\  N/D, \tau_b\  D_1).$$
\end{lem}

Let $\ceiling{t^m}$ be the map from $K[[t]]$ to $K[t]$ given by
replacing  $t^n$ with $0$ for all $n\ge m$. More precisely,
$$\ceiling{t^m} \sum_{n\ge 0} a_n t^n =\sum_{n=0}^{m-1} a_n t^n.$$
where $a_i\in K$ for all $i$. The following properties can be
easily checked for all $f,g\in K[[t]]$.
\begin{enumerate}
\item  $\ceiling{t^m} (f+g) =\ceiling{t^m} f +\ceiling{t^m} g $.
\item  $\ceiling{t^m} (fg) =\ceiling{t^m} (\ceiling{t^m} f \ceiling{t^m} g)$.
\item  If $0<k< m$ then $\ceiling{t^m} t^k f =t^k \ceiling{t^{m-k}} f$.
\item  If $g(0)\ne 0$, then $\ceiling{t^m} f/g =\ceiling{t^m} (\ceiling{t^m} f/ \ceiling{t^m} g)$.
\end{enumerate}


The main formula for our algorithm is the following, which is a
consequence of Theorem \ref{t-ppfraction-main1}. But we would like
to prove this result by using Laurent series expansion.
\begin{thm}\label{t-ppfraction-main-0}
Let $N,D\in K[t]$ and let $D=t^{m}E$ with $E\in K[t]$ and $E(0)\ne
0$. Then
$$ t^{m} \frr(N/D,t^{m})=\ceiling{t^{m}} \frac{N(t)}{E(t)}.$$
\end{thm}

\begin{pf}
Since $E(0)\ne 0$, $t^{m}$ and $E$ are relatively prime. Let
\begin{align}\label{e-ppfraction-th0}
\frac{N(t)}{D(t)}= p(t) +\frac{r(t)}{t^m} +\frac{s(t)}{E(t)}
\end{align}
be the ppfraction of $N/D$ with respect to $(t^m,E)$. Thus $\deg
(r(t))<m$, and $r(t)=t^m \frr(N/D,t^m)$.

Because $K(t)$ can be embedded into the field of Laurent series
$K((t))$, equation \eqref{e-ppfraction-th0} is also true as an
identity in $K((t))$.
 On the right-hand side of
equation \eqref{e-ppfraction-th0}, when expanded as Laurent series
in $K((t))$, the the second term contains only negative powers in
$t$, and the other terms contain only nonnegative powers in $t$.
Therefore, $r(t)/t^m$ equals the negative part of $N/D$ when
expanded as a Laurent series. More precisely, for $i=1,\dots , m$,
we have
$$[t^{-i}] \frac{N(t)}{D(t)} = [t^{-i}] \frac{ r(t)}{ t^m}.$$
This is equivalent to $[t^{m-i}] N(t)/E(t)= [t^{m-i}] r(t)$ for
$i=1,\dots, m$. Now $r(t)$ is a polynomial of degree at most
$m-1$, and $N(t)/E(t) \in K[[t]]$, so
$$r(t)= \ceiling{t^m} \frac{N(t)}{E(t)}.$$
\end{pf}

\begin{rem}
The idea of using Laurent series expansion to obtain part of the
partial fraction expansion appeared in the proof of \cite[Theorem
4.4]{gessel2}.
\end{rem}

Gessel observed that this same idea can also be used to compute
the polynomial part of a rational function, and that it is fast
when the polynomial part has small degree.
\begin{prop}\label{p-ppfraction-ira}
If $R(t)$ is a rational function in $K(t)$, then the polynomial
part $P(t)$ can be computed by the following equation.
$$t^{-1}P(t^{-1})= \frr( t^{-1} R(t^{-1}),t=0).$$
\end{prop}
\begin{pf}
Let $R(t)=P(t)+N(t)/D(t)$ be the ppfraction of $R(t)$, and let
$p=\deg (P)$, $d=\deg(D)$, and $n=\deg(N)$. Then $n<d$. Now we
have
$$t^{-1}R(t^{-1})=t^{-1}P(t^{-1})+t^{-1} N(t^{-1})/D(t^{-1}) =t^{-1}P(t^{-1})+
t^{d-n-1} \tilde{N}(t)/(\tilde{D}(t)),$$ where
$\tilde{D}(t)=t^dD(t^{-1})$, and similarly for $\tilde{N}(t)$.

Apply ppfraction expansion to the second term. Since
$\tilde{D}(t)$ has nonzero constant term, it is relatively prime
to $t^{p+1}$. Now it is clear that $t^{-1}P(t^{-1})$ is the
fractional part of $t^{-1}R(t^{-1})$ with respect to $t^{p+1}$.
\end{pf}

\begin{exa}
It is easy to check that
 $$R(t)={\frac {{t}^{3}+2\,{t}^{2}-3\,t+4}{{t}^{2}-4\,t+2}}
=t+6+{\frac {-8+19\,t}{{t}^{2}-4\,t+2}}.$$ Now we compute the
polynomial part of $R(t)$ by Proposition \ref{p-ppfraction-ira}.
\end{exa}
\begin{align*}
t^{-1}R(t^{-1})& ={\frac {1+2\,t-3\,{t}^{2}+4\,{t}^{3}}{t^2\left
(1-4\,t+2\,{t}^{2}\right
)}} \\
t^2 \frr( t^{-1} R(t^{-1}),t^2) &= \ceiling{t^2} {\frac
{1+2\,t-3\,{t}^{2}+4\,{t}^{3}}{\left (1-4\,t+2\,{t}^{2}\right
)}}\\
&=\ceiling{t^2} \frac{1+2t}{1-4t} =1+6t.
\end{align*}
So we obtain that the polynomial part of $R(t)$ is $t+6$.

\vspace{3mm}

Note that when expanded as Laurent series in $t$, we have
$$\ceiling{t^{m_0}} \frac1{(t-a_i)^{m_i}}=\sum_{j=0}^{m_0-1} (-1)^{m_i}
\binom{m_i-1+j}{j} \frac{t^j}{a_i^{m_i+j}}.$$ Hence by Theorem
\ref{t-ppfraction-main-0}, we get
\begin{cor}\label{c-ppfraction-0}
Let $N\in K[t]$ and $D=t^{m_0}(t-a_1)^{m_1}\cdots (t-a_k)^{m_k}$
with all the $a_i$'s distinct and not equal to $0$. Then
$$t^{m_0} \frr\left(\frac{N}{D},t^{m_0}\right)=\ceiling{t^{m_0}}
Ns_1\cdots s_k,$$ where
$$s_i=\sum_{j=0}^{m_0-1} (-1)^{m_i}
\binom{m_i-1+j}{j} \frac{t^j}{a_i^{m_i+j}}.$$
\end{cor}

Therefore, combining Theorem \ref{t-ppfraction-main1}, Lemma
\ref{l-ppfraction-tr} and Corollary \ref{c-ppfraction-0}, we
obtain an algorithm for computing the partial fraction
decomposition of a proper rational function of the general form
$$F(t)=\frac{N(t)}{(t-a_1)^{m_1}\cdots (t-a_k)^{m_k}}.$$

\begin{enumerate}
\item Let $S:=0$
\item For $i$ from $1$ to $k$ do $G(t) :=F(t+a_i),$
 $S:=S+\tau_{-a_i}\frr (G(t),t^{m_i} )$ next $i$.
 \item Return $S$.
\end{enumerate}

It was stated in \cite[Ch. 2.4]{xinthesis} that the computation of
$\frr(G(t),t^{m_i})$ for all $i$ will take time
$O(k(m_1^{1.58}+m_2^{1.58}+\cdots +m_k^{1.58}))$. However, this
estimate does not show $O(M^2)$ time. For example, $M=2k$, and
$m_1=m_2=\cdots =m_{k-1}$, and $m_k=M-k+1$.

\begin{thm}\label{t-algorithm0}
Let $M$ be the degree of the denominator of a rational function.
The above algorithm for partial fraction decomposition can be
executed in $O(M^2)$ time.
\end{thm}

The proof of this theorem, which will be given later, uses the
fact that manipulations in $K[t]/\langle t^m\rangle $ are fast.
Now let us estimate the computational time of manipulations in the
quotient ring $K[t]/\langle t^m\rangle $.

The following is a well-known result by the method of divide and
conquer. See, e.g., \cite[Property 36.1]{sedgewick}.
\begin{prop}
Let $R(m)$ be the time for computing the product of two
polynomials of degree less than $m$. Then $R(m)=O(m^{1.59})$.
\end{prop}

\begin{rem}
In the proof of Theorem \ref{t-algorithm0}, we only need the
obvious upper bound $R(m)=O( m^2)$. The above proposition shows
that our algorithm can be accelerated.
\end{rem}

Since most of our estimations use the method of divide and
conquer, it is better to introduce it here. In what follows, we
shall always assume that $m$ is a power of $2$ for simplicity. The
estimation of $R(m)$ follows from the following observation.
Bisect $P(t)$ as $P(t)= P_1(t)+t^{m/2}P_2(t)$, and bisect $Q(t)$
as $Q(t)=Q_1(t)+t^{m/2} Q_2(t)$. Then
\begin{align}
\label{ceilingn}
PQ=P_1Q_1+t^{m/2}((P_1+Q_1)(P_2+Q_2)-P_1Q_1-P_2Q_2)+t^{m}P_2Q_2.
\end{align}
which shows that we need only three polynomial multiplications.
This gives the recurrence $R(m)=3R(m/2)$ and that
$R(m)=O(m^{\log_2 3})=O(m^{1.59})$.

\begin{lem}\label{l-0-multiply} Let $P(t)$ and $Q(t)$ be two polynomials of
degree $m-1$. Then $\ceiling{t^m}  P(t)Q(t)$ can be computed in no
more than $R(m)$ time.
\end{lem}

The proof of this lemma is trivial.

%

\begin{lem}
The computation of $\ceiling{t^m}P/Q$, where $Q(0)\ne 0$, takes no
more than $2R(m)$ time.
\end{lem}
\begin{pf} We use the method of divide and conquer. Let $T(m)$ be the time
for the computation in question.

Suppose that $\ceiling{t^m}P/Q=Z$. Bisect $P,Q,Z$ as
$P=P_1+t^{m/2}P_2$, $Q=Q_1+t^{m/2}Q_2$, and $Z=Z_1+t^{m/2}Z_2$.
Then
$$\ceiling{t^m}\frac{P_1+t^{m/2}P_2}{Q_1+t^{m/2}Q_2}=Z_1+t^{m/2}Z_2.$$

Obviously $\ceiling{t^{m/2}} P_1/Q_1=Z_1$. So it will take
$T(m/2)$ time to find $Z_1$. To find $Z_2$, we use the formula
(from direct algebraic computation).
$$Z_2=\ceiling{t^{m/2}}\frac{(P_1-Q_1Z_1)/t^{m/2}+P_2-Q_2Z_1}{Q_1}.$$
Therefore, we get the recurrence $T(m)=T(m/2)+2R(m/2)+T(m/2)$,
where the first summand is for $Z_1$, and the rest is for $Z_2$.
Using this recurrence, it is easy to see that $T(m)$ is no more
than $2R(m)$.
\end{pf}

To apply the above lemma, we need the expanded representation of
$P(t)$ and $Q(t)$.

\begin{lem}
Suppose that $Q(t)$ is the product of $m$ linear factors. Then it
takes no more than $R(m)$ time to expand $Q(t)$.
\end{lem}
\begin{pf}
Let $U(m)$ denote the time for expanding the products of $m$
linear factors. Factor $Q(t)$ as $Q(t)=Q_1(t)Q_2(t)$, where
$Q_1(t)$ consists of the first $m/2$ factors. Then it will take
$U(m/2)$ time to expand $Q_1(t)$, and $U(m/2)$ time to expand
$Q_2(t)$, and then $R(m/2)$ time to get the final expansion. Thus
$U(m)=2U(m/2)+R(m/2)$. This implies that $U(m)$ is approximately
equal to $R(m)$.
\end{pf}

We might be able to speed up the expansion in the above lemma by
the following lemma.

\begin{lem}
The expansion of $(t-a)^m$ takes $O(m)$ time.
\end{lem}
\begin{pf}
  This lemma follows from the binomial theorem
  $$(t-a)^m=\sum_{i=0}^{m}t^i(-a)^{m-i}\binom{m}{i},$$
  and the fact that the ratios of consecutive summands are simple.
\end{pf}

Using the binomial theorem, it is easy to see the following.
\begin{lem}
Suppose the degree of $N(t)$ is less than $M$. Then the expansion
of $\ceiling{t^m} N(t+a)$ takes no more than $O(m\cdot \text{the
number of nonzero terms in $N(t)$})$, which is no more than
$O(Mm)$ time.
\end{lem}

Now we estimate the computational time for $\frr(N/D,m_0)$ in
Theorem \ref{t-ppfraction-main-0}.
\begin{prop}
Let $M=m_0+m_1+\cdots +m_k$. Then it takes $O(Mm_0^{0.59})$ time
to compute $\ceiling{t^{m_0}} \frac{N(t)}{(t-a_1)^{m_1}\cdots
(t-a_k)^{m_k}}$.
\end{prop}
\begin{pf}
We first expand the denominator by grouping every $(m_0-1)$
factors together. So we have about $(M-m_0)/(m_0-1)$ groups. It
takes $R(m_0)$ time for expanding the products for each group, and
then about $(M-m_0)/(m_0-1)$ multiplications when taking
$\ceiling{t^{m_0}}$. Thus the total time for this expansion is
about $2 M/m_0 R(m_0)$.

Denote by $E(t)$ the resulting expansion. Now it will take about
$2R(m_0)$ time to compute $\ceiling{t^{m_0}} \frac{N(t)}{E(t)}$.
Therefore the total time for the final answer is \\
$O(M/m_0)R(m_0)=O(Mm_0^{0.58})$.
\end{pf}

\begin{pf*}{Proof of Theorem \ref{t-algorithm0}.} It will take $O(Mm_i)$ time
for finding the expansion of $N(t+a_i)$, and will take
$O(Mm_i^{0.58})$ time for finding
$\ceiling{t^{m_i}}t^{m_i}N(t+a_i)/D(t+a_i)$. Therefore, the total
time for finding $\frr(N/D,x=a_i)$ takes $O(Mm_i)$ time. Summing
on all $i$, we see that it takes $O(M^2)$ time to find the partial
fraction decomposition of $N/D$.
\end{pf*}

This new algorithm also enables us to work with some
 difficult rational functions by hand.

\begin{exa}
Compute the partial fraction expansion of $f(t)$, where
$$f(t)= \frac{t}{(t+1)^2(t-1)^3(t-2)^5}.$$
\end{exa}

\begin{pf*}{Solution.} Clearly, the polynomial part of $f(t)$ is $0$. Although
applying Corollary \ref{c-ppfraction-0} is faster, we compute the
fractional part of $f(t)$ at $t=-1$ and $t=1$ differently. For the
fractional part of $f(t)$ at $t=-1$, we apply $\tau_{-1}$, and
compute $\frr (f(t-1),t^2)$ by Theorem \ref{t-ppfraction-main-0}.
We have
\begin{align*}
t^2 \frr(f(t-1),t^2) &= \ceiling{t^2} \frac{t-1}{(t-2)^3(t-3)^5}\\
  &= \ceiling{t^2} \frac{ t-1}{ (-8+12t)((-3)^5+3^4\cdot 5t)} \\
  &= \ceiling{t^2} \frac{t-1}{8\cdot 3^5(1-19/6t)} \\
  &= \ceiling{t^2} \frac{(t-1)(1+19/6t)}{8\cdot 3^5} = -\frac{1}{8\cdot 3^5}(1+\frac{13t}{6}).
\end{align*}
Thus
$$\frr(f(t), (t+1)^2)= -\frac{1}{2^3\cdot 3^5 (t+1)^2}-\frac{13}{2^4 \cdot 3^6 (t+1)}.$$

Similarly, we can compute the fractional part of $f(t)$ at $t=1$.
We have
\begin{align*}
t^3 \frr(f(t+1),t^3) &= \ceiling{t^3} \frac{t+1}{(t+2)^2(t-1)^5}\\
  &= \ceiling{t^3} \frac{ t+1}{ (t^2+t+4)(-10t^2+5t-1)} \\
  &= \ceiling{t^3} \frac{t+1}{-4+16t-21t^2} \\
  &= -\frac{1}{4}\ceiling{t^3} (t+1)(1+4t-\frac{21}{4}t^2 +16t^2) \\
  &=-\frac{1}{4}(1+5t+\frac{59}{4} t^2).
\end{align*}
Thus
$$\frr(f(t), (t-1)^3)=  -\frac{1}{4(t-1)^3}-\frac{5}{4(t-1)^2}-
\frac{59}{16(t-1)}.$$ The fractional part of $f(t)$ at $t=2$ can
be obtained similarly, but it is better to use Corollary
\ref{c-ppfraction-0}. In fact, this computation becomes quite
complicated. Although it is still possible to work by hand, we did
use Maple.
\begin{align*}
&t^5\frr(f(t+2),t^5) \\&=\ceiling{t^5} (t+2)
\left(\frac{1}{9}-{\frac {2t}{27}}+\frac{{t}^{2}}{27}-{\frac
{4{t}^{3}}{243}}+{\frac {5{t}^{4}}{729}} \right)
\left(1-3t+6{t}^{2}-10{t}^{3}+15{t}^{4}
 \right)
 \\&= \frac{2}{9}-{\frac
{19}{27}}t+{\frac {13}{9}}{t}^{2}-{\frac {593}{243}}{
t}^{3}+{\frac {2689}{729}}{t}^{4}.
\end{align*}
Applying Theorem \ref{t-ppfraction-main1}, we get the partial
fraction expansion of $f(t)$, which is too lengthy to be worth
giving here.
\end{pf*}

\section{Algorithm for a General Field and Full Partial Fraction Decompositions}

When $K$ is a general field, e.g., the field of rational numbers
$\QQ$, linear transformations will not help. Manipulations in
$K[t]/\langle D_1(t)\rangle $ are not as good as the case of
$D_1(t)=x^m$. But we still have an $O(M^2)$ algorithm.

\begin{prop}\label{p-ppfraction-k}
Suppose that $D_1,\dots, D_k\in K[t]$ are pairwise relatively
prime, and $D=D_1\cdots D_k$. If $\deg(N)<\deg(D)$, then the
ppfraction decomposition of $N/D$ with respect to $D_1,\dots ,D_k$
can be computed in $O(\deg(D)^2)$ time.
\end{prop}

The proof of this proposition will be given later. Now assume that
$D_i=p_i^{a_i}$, and $\deg(D_i)=m_i$. It is easy to show that the
partial fraction decomposition of $r_i/D_i$, where
$\deg(r_i)<m_i$, can be computed in $O(m_i^2)=O(Mm_i)$ time. Thus
the above argument and Proposition \ref{p-ppfraction-k} will give
us the following.

\begin{thm}
Suppose that $\deg(N)<\deg(D)$, and we are given a factorization
$D=p_1^{a_1}\cdots p_k^{a_k}$ of $D$ into primes in $K[t]$. Then
the partial fraction decomposition of $N/D$ takes $O(\deg(D)^2)$
time.
\end{thm}

In order to prove Proposition \ref{p-ppfraction-k}, we need to
estimate manipulations in $K[t]/\xmod{D}$ for a given polynomial
$D$. In most situations, we need the unique representative of
$N+\xmod{D}$ that has degree less than $\deg{D}$. We denote by
$\ceiling{D} N$ this representative, which is also known as the
remainder of $N$ when divided by $D$.

The following estimations are obvious. The computational time
refers to the number of multiplications of two elements in $K$.
Time spent on additions is omitted.

\begin{enumerate}
    \item The computational time for expanding $PQ$ for any two
polynomials $P$ and $Q$ is no more than $(\deg(P)+1)(\deg(Q)+1)$
time.

    \item For two given polynomials $N$ and $D$, the division algorithm for
finding $p$ and $r$ such that $N=pD+r$ with $\deg(r)<\deg(D)$
takes  no more than $(\deg(N)-\deg(D))\cdot (\deg(D)+1) $ time.

    \item Suppose that $P$ and $Q$ are two polynomials of degree less
than $\deg(D)$. Then the computation of $\ceiling{D} PQ$ takes no
more than $2\deg(D)^2 $ time.
\end{enumerate}

\begin{lem}\label{l-euclidean}
Suppose $P$ is relatively prime to $D$ and $\deg(P)<\deg(D)$. Then
the computation of $\ceiling{D} 1/P$ takes $O(\deg(D)^2)$ time.
\end{lem}
This estimation is obtained by the extended Euclidean algorithm
for polynomials. See, e.g., \cite{moenck}, which says that an
$O(\deg(D)\log^r(D))$ algorithm exists, where $r$ is a fixed
number.

\begin{lem}
  Suppose that $\deg(D_1)=m$ and $\deg(D_2\cdots D_k)=M$. Then the
  computation of $\ceiling{D_1} D_2\cdots D_k$ takes no more than
  $ 4(M+m)m$ time.
\end{lem}
\begin{pf}
Denote by $V(M)$ the computational time described in the lemma. We
shall prove that $V(M)\le \max \{\, 4Mm-2m^2,Mm+m^2\, \}$, which
implies the lemma. The proof is in two parts. The first part deals
with the case when $M\le 2m$, and the second part deals with the
case when $M\ge m$. Note that there is an overlap.

We first show that the expansion of $D_2\cdots D_k$ takes no more
than $M^2/2$ time by induction on $M$. This claim is clearly true
for small $M$, e.g., $M=1,2$. Now suppose the claim is true for
all $l\le M$. Then the expansion of $D_2\cdots D_k$ can be
obtained by first expanding $D_2\cdots D_{k-1}$ (of degree $M_1$),
then multiplying it by $D_k$ (of degree $M_2$). The computational
time is (by induction) no more than $M_1^2/2+M_1M_2 \le
(M_1+M_2)^2/2=M^2/2$. Thus for $M\le 2m$, we get $V(M)\le
M^2/2+m^2\le Mm+m^2$ by expanding $D_2\cdots D_k$, and then taking
$\ceiling{D_1}$.

We claim that for all $M\ge m$, $V(M)\le 4Mm-2m^2$, and prove the
claim by induction on $M$. The claim follows from the inequality
$M^2/2+m^2< 4Mm-2m^2$ when $M\le 7m$ by the first part. For $M\ge
7m$, we can separate $D_2\cdots D_k$ into two products of degree
$M_1$ and $M_2$ respectively. We can assume that $M_1\ge m$ and
$M_2\ge m$, for otherwise, the degree of one of $D_i$ is larger
than $5m$, in which case the claim is easily seen to be true. Now
we compute the remainder of each product, and then compute the
resulting product and compute the remainder. This process takes
time
$$V(M_1)+V(M_2)+2m^2\le 4M_1m-2m^2+4M_2m-2m^2+2m^2=4Mm-2m^2.$$
This completes the proof.
\end{pf}

\begin{prop}\label{p-remainder-1}
Suppose that $\deg(D)=M$, $\deg(N)\le M$, $\deg(D_1)=m$, and that
$D=D_1\cdots D_k$ is a factorization of $D$ into relatively prime
factors. Then the computation of $\ceiling{D_1} N/(D_2\cdots D_k)$
takes $O(Mm)$ time.
\end{prop}
\begin{pf}
We first compute $\ceiling{D_1} D_2\cdots D_k$, and denote the
result by $D'$. This step takes no more than $4Mm$ time. Then we
compute $\ceiling{D_1} N$, and denote the result by $N'$. This
step takes no more than $Mm$ time. Finally we compute
$\ceiling{D_1} N'/D'$. This step takes $O(m^2)=O(Mm)$ time by
Lemma \ref{l-euclidean}. So the total time is
$4Mm+Mm+O(Mm)=O(Mm)$.
\end{pf}

\begin{pf*}{Proof of Proposition \ref{p-ppfraction-k}.} By Theorem
\ref{t-ppfraction-main1}, the numerator of $\frr(N/D,D_i)$ is
given by
$$r_i=\ceiling{D_i} ND_i/(D_1D_2\cdots
D_k).$$ The computation of $r_i$ takes $O(M\deg(D_i))$ time by
Proposition \ref{p-remainder-1}. Summing on all $i$ we get the
total computational time for the ppfraction of $N/D$ with respect
to $D_1,\dots ,D_k$, which is $\sum_{i=1}^kO(Mm_i)=O(M^2)$.
\end{pf*}

\begin{exa}
Compute the fractional part of $f(t)$ with respect to $t^2-t+2$,
where
$$f(t)=\frac{t^2}{(t^2-2t-1)^2(t^2-t+2)}.$$
\end{exa}
\begin{pf*}{Solution.} Let $p(t)=t^2-t+2$. Then we need to compute
$\ceiling{p(t)} t^2/(t^2-t-1)^2$. In the following computation, we
shall always replace $t^2$ with $t-2$.
\begin{align*}
\ceiling{p(t)} \frac{t^2}{(t^2-2t-1)^2} &= \ceiling{p(t)}
\frac{t-2}{(-t-3)^2}=\ceiling{p(t)}\frac{t-2}{7(t+1)}\\
&=\ceiling{p(t)}\frac{(t-2)^2}{-28}=\frac{-3t+2}{-28},
\end{align*}
where we used the fact that  $(t+1)(t-2)=t^2-t-2=p(t)-4$.
Therefore
$$\frr({f(t),p(t)})=\frac{3t-2}{-28(t^2-t+2)}.$$
\end{pf*}

 In Maple, the full partial fraction expansion of a
rational function will involve a form like
$$\sum_{\alpha= \text{\rm root of } p(t)}\sum_{j=1}^m \frac{h_j(\alpha)}{(t-\alpha)^j},$$
where $p(t)$ is a prime polynomial, and $h_j(t)$ will be a
polynomial of degree less than $\deg(p(t))$. This expansion is
useful in some situations. We can get this kind of expansion by
applying Theorem \ref{t-ppfraction-main-0}. This is best
illustrated by an example.

\begin{exa}
Compute the full partial fraction expansion of $f(t)$, where
$$f(t)=\frac{t}{(t^2-t-1)^2(t^2-t+2)}.$$
\end{exa}

\begin{pf*}{Solution.} Suppose that $\alpha$ is a root of the prime polynomial
$p(t):=t^2-t-1$. Since $K(\alpha)$ is a field, and
$\alpha^2=\alpha+1$, we can use this relation to get rid of all
terms containing $\alpha^n$ for $n\ge 2$. Because $p(t)$ is a
prime polynomial, $\alpha$ can only be a simple root of $p(t)$.
Then $t$ divides $p(t+\alpha)$ and $p(t+\alpha)/t$ has nonzero
constant term. In the present example,
$$p(t+\alpha)= (t+\alpha)^2-(t+\alpha)-1= t(t+2\alpha-1).$$
Note that the constant term of $p(t+\alpha)$ is always $0$.

Clearly, $\tau_\alpha (t^2-t+2)$ has nonzero constant term, for
otherwise $t^2-t+2$ will not be relatively prime to $p(t)$. In the
present situation,
$$ (t+\alpha)^2-(t+\alpha)+2=t^2+(2\alpha -1)t+3.$$
By Lemma \ref{l-ppfraction-tr} and Theorem
\ref{t-ppfraction-main-0}, we can work in $K(\alpha)[[t]]$.
\begin{align*}
\ceiling{t^2} t^2 f(t+\alpha) &= \ceiling{t^2}
\frac{t+\alpha}{(t+2\alpha -1)^2(t^2+(2\alpha-1)t+3)}
\\
&= \frac{1}{15} \ceiling{t^2} \frac{t+\alpha}{1+(2\alpha-1)11t/15}\\
&= \frac{1}{15} \ceiling{t^2} (t+\alpha)(1-11(2\alpha-1)t/15) \\
&= \frac{1}{15}\alpha +\frac{(-11\alpha+7)}{15^2}t.
\end{align*}
Thus the fractional part of $f(t)$ at $\alpha$ that satisfies
$p(\alpha)=0$ can be written as
$$  \frac{\alpha}{15(t-\alpha)^2}
+\frac{(7-11\alpha)}{225(t-\alpha)}.$$ Similarly, the fractional
part of $f(t)$ at $\beta$ that satisfies $\beta^2-\beta+2=0$ can
be written as
$$ \left({\frac {4}{63}}-{\frac {1}{63}}\beta\right)\left (t-\beta
\right )^{-1}
 .$$
Together with the fact that the polynomial part of $f(t)$ is
clearly $0$, the full partial fraction expansion of $f(t) $ is
hence
$$f(t)= \sum_{\alpha^2-\alpha-1=0} \left[\frac{\alpha}{15(t-\alpha)^2}
+\frac{(7-11\alpha)}{225(t-\alpha)}\right]
+\sum_{\beta^2-\beta+2=0} \frac{4-\beta}{63(t-\beta) }.$$
\end{pf*}

Of course we can first compute the partial fraction decomposition
of $f(t)$ and then compute its full partial fraction
decomposition.

\vspace{3mm} {\bf Acknowledgment.} The author is very grateful to
his advisor Ira Gessel.

\bibliographystyle{elsart-harv}


\begin{thebibliography}{5}
\expandafter\ifx\csname
natexlab\endcsname\relax\def\natexlab#1{#1}\fi
\expandafter\ifx\csname url\endcsname\relax
  \def\url#1{\texttt{#1}}\fi
\expandafter\ifx\csname
urlprefix\endcsname\relax\def\urlprefix{URL }\fi

\bibitem[{Gessel(1997)}]{gessel2}
Gessel, I.~M., 1997. Generating functions and generalized
{D}edekind sums.
  Elec. J. Comb. 4~(2), Wilf Festschrift, R11.

\bibitem[{Kovacic(1986)}]{kovacic}
Kovacic, J.~J., 1986. An algorithm for solving second order linear
homogeneous
  differential equations. J. of Symbolic Computation 13, 3--43.

\bibitem[{Moenck(1973)}]{moenck}
Moenck, R.~T., 1973. Fast computation of gcds. In: Proceedings of
the fifth
  annual ACM symposium on Theory of computing. Austin, Texas, United States,
  pp. 142--151.

\bibitem[{Sedgewick(1988)}]{sedgewick}
Sedgewick, R., 1988. Algorithms, 2nd Edition. Addison-Wesley, New
York.

\bibitem[{Xin(2004)}]{xinthesis}
Xin, G., 2004. {T}he {R}ing of {M}alcev-{N}eumann {S}eries and
{T}he {R}esidue
  {T}heorem. Ph.D. thesis, Brandeis University.

\end{thebibliography}

\end{document}